\documentclass[11pt]{article}
\newtheorem{Lemma}{Lemma}
\newtheorem{Proposition}[Lemma]{Proposition}

\newtheorem{Theorem}[Lemma]{Theorem}

\newcommand{\cd}{\ \stackrel{d}{\rightarrow} \ }
\newcommand{\cp}{\ \stackrel{p}{\rightarrow} \ }
\newcommand{\TT}{\mbox{${\cal T}$}}

\newcommand{\FF}{\mbox{${\cal F}$}}

\newcommand{\eps}{\varepsilon}

\newcommand{\sfrac}[2]{{\textstyle\frac{#1}{#2}}}

\usepackage{amssymb}
\newcommand{\len}{{\rm len}}
\newcommand{\perc}{{\rm perc}}
\newcommand{\exc}{{\rm exc}}

\newcommand{\Rbold}{{\mathbb{R}}}
\newcommand{\Pb}{{\mathbb{P}}}
\newcommand{\Eb}{{\mathbb{E}}}

\newcommand{\Zbold}{{\mathbb{Z}}}
\newcommand{\Cbold}{{\mathbb{C}}}

\def\ut{\underline{t}}

\def\ind{{\rm 1\hspace{-0.90ex}1}}
\def\utx{\underline{t-x}}

\usepackage{epsfig}
\author{David J. Aldous\thanks{Research supported by N.S.F. Grant
DMS0203062}
\\
\\
\\
       University of California\\
       Department of Statistics\\
        367 Evans Hall \# 3860\\
       Berkeley CA 94720-3860\\
       aldous@stat.berkeley.edu
\and Charles Bordenave\thanks{Research supported in part by the EuroNGI network}\\
\\
\\
\'Ecole Normale Sup\'erieure\\
D\'epartement d'Informatique \\
45 rue d'Ulm\\
75230 Paris cedex 5 \\
charles.bordenave@ens.fr
\and Marc Lelarge\thanks{Research supported
by Science Foundation Ireland Grant SFI04/RP1/I512}\\
\\
\\
       INRIA-ENS\\
       45 rue d'Ulm\\
       75230 Paris Cedex 5\\
       marc.lelarge@ens.fr}

\title{Near-Minimal Spanning Trees: a Scaling Exponent in Probability Models}

\begin{document}

\maketitle

\begin{abstract}
We study the relation between the minimal spanning tree (MST) on
many random points and the ``near-minimal" tree which is optimal
subject to the constraint that a proportion $\delta$ of its edges
must be different from those of the MST.  Heuristics suggest that,
regardless of details of the probability model, the ratio of lengths
should scale as $1 + \Theta(\delta^2)$.  We prove this {\em scaling
result} in the model of the lattice with random edge-lengths and in
the Euclidean model.
\end{abstract}


\noindent
{\bf Mathematical subject codes:}
05C80;
60K35;
68W40

\section{Introduction}
This paper gives details of one aspect of the following broad project
\cite{me103}.
Freshman calculus tells us how to find a minimum $x_*$ of a smooth
function $f(x)$: set the derivative $f^\prime(x_*) = 0$
and check $f^{\prime \prime}(x_*) > 0$.
The related series expansion tells us, for points $x$ near to $x_*$,
how the distance $\delta = |x-x_*|$ relates to the difference
$\eps = f(x) - f(x_*)$ in $f$-values:
$\eps$ scales as $\delta^2$.  This {\em scaling exponent} $2$
persists for functions 
$f: \Rbold^d \to \Rbold$:
if $x_*$ is a local minimum and $\eps(\delta) := \min\{f(x)-f(x_*):|x-x_*| =
\delta)$,
then
$\eps(\delta)$ scales as $\delta^2$
for a generic smooth function $f$.

Combinatorial optimization, exemplified by the
{\em traveling salesman problem}
(TSP),
is traditionally viewed as a quite distinct subject, with theoretical
analysis focussed on the number of steps that algorithms require to find
the optimal solution.
To make a connection with calculus,
compare an arbitrary tour ${\bf x}$ through $n$ points with the optimal (minimum-length)
tour ${\bf x}_*$ by considering the two quantities
\begin{eqnarray*}
\delta_n({\bf x}) = \mbox{\{number of edges in {\bf x}
but not in {\bf x}$_*$\}}/n \\
\eps_n({\bf x}) = \mbox{\{length difference between {\bf x}
and {\bf x}$_*$\}}/s(n)
\end{eqnarray*}
where $s(n)$ is the length of the minimum length tour.
Now define $\eps_n(\delta)$ to be the minimum value of
$\eps_n({\bf x})$ over all tours ${\bf x}$
for which $\delta_n({\bf x}) \geq \delta$.
Although the function $\eps_n(\delta)$ will depend on $n$ and
the problem instance, we anticipate that for typical instances
drawn from a suitable probability model
it will converge in the
$n \to \infty$ limit to some deterministic function
$\eps(\delta)$.  The {\em universality\/}
paradigm from statistical physics \cite{kad00} suggests there
might be a scaling exponent $\alpha$ defined by
\[ \eps(\delta) \sim \delta^\alpha \mbox{ as } \delta \to 0 \]
and that the exponent should be robust under model details.

There is fairly strong evidence \cite{me103} that for TSP
the scaling exponent is $3$.
This is based on analytic methods in a {\em mean-field}
model of interpoint distances
(distances between pairs of points are random, independent for
different pairs, thus ignoring geometric constraints)
and on Monte Carlo simulations for random points in $2$, $3$ and $4$
dimensional space.
The analytic results build upon a recent probabilistic
reinterpretation \cite{me94}
of work of Krauth and M{\'e}zard
\cite{KM87}
establishing the average length of mean-field TSP tours.
But neither part of these TSP assertions is rigorous, and indeed rigorous proofs
in $d$ dimensions seem
far out of reach of current methodology.
In contrast, for the
{\em minimum spanning tree} (MST) problem, a standard algorithmically easy
problem, a simple heuristic argument (section \ref{sec-heuristic})
strongly suggests that the scaling exponent is $2$ for any reasonable
probability model.
The goal of this paper is to work through the details of a rigorous proof.

Why study such scaling exponents?
For a combinatorial optimization problem,
a larger exponent means that there are more near-optimal solutions,
suggesting that the algorithmic problem of finding the optimal
solution is intrinsically harder.
So scaling exponents may serve to separate
combinatorial optimization problems of an appropriate type
into a small set of classes of increasing difficulty.
For instance, the {\em minimum matching} and {\em minimum Steiner tree} problems
are expected to have scaling exponent $3$,
and thus be in the same class as TSP in a
{\em quantitative} way, as distinct from their qualitative similarity
as NP-complete problems under worst-case inputs.
In contrast, algorithmically easy problems are expected to have
scaling exponent $2$, analogously to the ``calculus" scaling exponent.
One plausible explanation is that the near-optimal solutions in such problems
differ from the optimal solution via only ``local changes", each local change affecting only a number
of edges which remains $O(1)$ as $\delta \to 0$.

\subsection{Background}
\label{sec-back}
Steele \cite{steele97}
and Yukich \cite{yukich-book}
give general background concerning combinatorial optimization
over random points.

A {\em network} is a graph whose edges $e$ have positive real
{\em lengths} $\len (e)$.
Let $G$ be a finite connected network.
Recall the notion of a {\em spanning tree} (ST) $T$ in $G$.
Identifying $T$ as a set of edges, write
$\len (T) = \sum_{e \in T} \len (e)$.
A {\em minimal spanning tree} (MST) is a ST of minimal length;
such a tree always exists but may not be unique.
The classical greedy algorithm ({\em Kruskal's algorithm} \cite{CL86}) for constructing a MST yields two
fundamental properties which we record without proof in Lemma \ref{L0}.

Let $G_t$ be the subnetwork consisting of those edges $e$ of $G$ with
$\len (e) < t$.
For arbitrary vertices $v,w$ define
\begin{equation} \perc (v,w) = \inf \{t: \mbox{ $v$ and $w$ in same component of $G_t$
} \}. \label{def-perc}
\end{equation}
For an edge $e = (v,w)$ of $G$ write $\perc(e) = \perc(v,w) \leq
\len(e)$ and also define the {\em excess}
\[ \exc (e) = \len(e) - \perc(e) \geq 0 .\]
\begin{Lemma}
\label{L0}
Suppose all the edge-lengths in $G$ are distinct.\\
(a) There is a unique MST, say $T$, and it is specified by the criterion
\[
e \in T \mbox{ if and only if } \exc(e) = 0 . \]
(b) For any vertices $v,w$
\[ \perc (v,w) = \max \{ \len(e): e
\mbox{ on path from $v$ to $w$ in $T$}
\} . \]
\end{Lemma}

\subsection{The heuristic argument}
\label{sec-heuristic}
Given a probability model for $n$ random points and their interpoint
lengths, define a measure $\mu_n(\cdot)$ on $ (0,\infty)$ in
terms of the expectation
\[
\mu_n(0,x) = \frac{1}{n} \Eb \left| \{\mbox{ edges }e : 0 < \len(e)
- \perc(e) < x\,\} \right|.
\]
For any reasonable model with suitable scaling of edge-lengths we expect an $n \to \infty$
limit measure $\mu(\cdot)$, with a density $f_{\mu}(x)=d\mu/dx$ having
a non-zero limit $f_{\mu}(0^+)$
as $x \downarrow 0$.


%
Now modify the MST by adding an edge $e$ with $\len(e) - \perc(e) =
b$, for some small $b$, to create a cycle; then delete the longest
edge $e^\prime \neq e$ of that cycle, which necessarily has
$\len(e^\prime) = \perc(e)$. This gives a spanning tree containing
exactly one edge not in the MST and having length greater by $b$.
Repeat this procedure with {\em every\/} edge $e$ for which
$0<\len(e) - \perc(e) < \beta$, for some small $\beta$. For large
$n$, the number of such edges should be $n\mu_n(0,\beta)\approx
n\,f_{\mu}(0^+)\beta$ to first order in $\beta$, and assuming there is
negligible overlap between cycles, each of the new edges will
increase the tree length by $\sim\beta/2$ on average.  
So we expect (Lemma \ref{LxiW})
\[ \delta(\beta) \sim f_{\mu}(0^+)\beta, \quad
\eps(\beta) \sim f_{\mu}(0^+)\beta^2/2 . \]
This construction should yield essentially the
minimum value of $\eps$ for given $\delta$, so we expect
\begin{equation}
\eps(\delta) \sim \frac{\delta^2}{2 f_{\mu}(0^+)}
\label{nu+}
\end{equation}
and in particular we expect
the scaling exponent to be $2$.

\subsection{Results}
Our goal is to formalize the argument above in the
context of the following two probability models for
$n$ random points.
Fix dimension $d \geq 2$
(the case $d=1$ is of course rather special).

\paragraph{Model 1}
{\em The disordered lattice}.
Start with the discrete $d$-dimensional cube
$\Cbold^d_m = [1,2,\ldots,m]^d$,
so there are $n = m^d$ vertices and there are $2d$ edges at each non-boundary
vertex.
Then take the edge-lengths to be i.i.d. random variables $\xi_e$,
whose common distribution $\xi$ has finite mean and
some bounded continuous density function $f_\xi(\cdot)$.

\paragraph{Model 2}
{\em Random Euclidean}.
Take the continuum $d$-dimensional cube $[0,n^{1/d}]^d$ of volume $n$.
Put down $n$ independent uniformly distributed random points in this cube.
Take the complete graph on these $n$ vertices, with Euclidean distance as
edge-lengths.

The results of this paper will remain valid in a slightly more
general framework than Model 2 in which points are put down
independently at random in the cube $[0,n^{1/d}]^d$ with common
density $f(n^{1/d}x)$ on $\Rbold^d$, with $f$ having support on
$[0,1]^d$ and being bounded away from zero. To avoid technicalities,
we restrict ourselves to the case $f$ constant.

Each model is set up so that nearest-neighbor distances are order $1$
and the MST $T_n$ has mean length of order $n$.
To formalize the ideas in the introduction we define the random variable
\begin{equation}
\eps_n(\delta) :=
\min \left\{
\frac{\len(T_n^\prime) - \len(T_n)}{n}: |T_n^\prime \setminus T_n| \geq
\delta n \right\}
 \label{eps-def}
\end{equation}
where the minimum is over spanning trees $T_n^\prime$
and where $T_n^\prime \setminus T_n$ is the set of edges
in $T_n^\prime$ but not in $T_n$.
\begin{Theorem}
\label{T1}
In either model, we have
\begin{eqnarray*}
{(a)}\quad \limsup_{\delta \downarrow 0} \delta^{-2}
\limsup_n \Eb \eps_n(\delta) < \infty , 
\end{eqnarray*} and,
\[(b) \quad \liminf_{\delta \downarrow 0} \delta^{-2}
\liminf_n \Eb \eps_n(\delta) > 0 . \]
\end{Theorem}

\paragraph{Structure of the paper}  In Section
\ref{sec:finite}, we do calculations  in the finite models: we prove
Theorem \ref{T1} for Model 1 and part (a) of the theorem for Model
2. In Section \ref{sec:infinite}, we introduce the limit infinite
random network (limit in the sense of {\em local weak convergence}
\cite{me101}) and its associated minimal spanning forest. We show
how results from continuum percolation theory allow us to show part
(b) of Theorem \ref{T1} for Model 2. 

\section{Proofs for the finite network}\label{sec:finite}
\subsection{The upper bound: Model 1 with $d=2$}
We first consider Model 1 with $d = 2$
and then consider the other cases.

The upper bound rests upon a simple construction of
near-minimal spanning trees, illustrated in Figure 1.

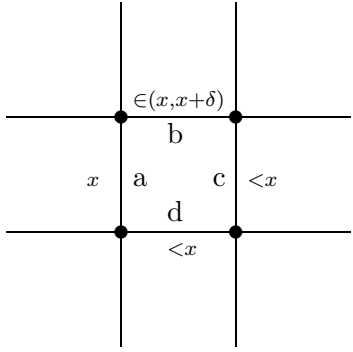
\begin{figure}
\begin{center}\setlength{\unitlength}{0.06in}
\begin{picture}(40,32)(0,0)
\put(10,0){\line(0,1){30}} \put(20,0){\line(0,1){30}}
\put(0,10){\line(1,0){30}} \put(0,20){\line(1,0){30}}
\put(10,10){\circle*{1.2}} \put(20,10){\circle*{1.2}}
\put(20,20){\circle*{1.2}} \put(10,20){\circle*{1.2}} \put(11,14){a}
\put(18,14){c} \put(14,17.8){b} \put(14,11){d}
\put(7,14){$\scriptstyle{x}$} \put(21,14){$\scriptstyle{<x}$}
\put(11,21){$\scriptstyle{\in (x,x+\delta)}$}
\put(14,8){$\scriptstyle{<x}$}
\end{picture}\end{center}
 \caption{A special configuration on
the $3 \times 3$ grid.}
\end{figure}

\noindent
The figure illustrates a particular kind of configuration.
There is a $4$-cycle of edges $abcd$ where, for some $x$,
\[ \len(a) = x, \ \len(b) \in (x, x+ \delta), \
\len(c)<x, \ \len(d) < x \] and where the eight other edges touching
the cycle have lengths $>x+\delta$. With such a configuration
(within a larger configuration on $\Cbold^2_m$), edges $adc$ are in
the MST, and edge $b$ is not. We can modify the minimal spanning
tree by removing edge $a$ and adding edge $b$; this creates a new
spanning tree whose extra length equals $\len(b) - x$.

Thus given a realization of the edge-lengths on the
$m \times m$ discrete square,
partition the square into adjacent $3 \times 3$ regions;
on each region where the configuration is as in Figure 1,
make the modification above.
This changes the MST $T_n$ into a certain near-minimal spanning
tree $T_n^\prime$.
On each $3 \times 3$ square,
the probability of seeing the Figure 1 configuration equals
\[ q(\delta) :=
\int_0^\infty f(x) (F(x+\delta)-F(x))F^2(x)(1-F(x+\delta))^8 \ dx
. \]
Here $f$ and $F$ are the density and distribution functions
of edge-lengths.
And the (unconditioned) increase in edge-length of spanning tree
caused by the possible modification equals
\[ r(\delta):=
\int_0^\infty f(x)
\left(\int_x^{x+\delta} (y-x)f(y)dy\right)
F^2(x)(1-F(x+\delta))^8 \ dx . \]
Letting $n \to \infty$ with fixed $\delta$,
and using the weak law of large numbers,
\begin{eqnarray}
 n^{-1} |T^\prime_n \setminus T_n| &\cp& \sfrac{1}{9} q(\delta) \label{con-1}\\
 n^{-1} (\len(T^\prime_n) - \len(T_n)) &\cp& \sfrac{1}{9} r(\delta) \label{con-2} .
\end{eqnarray}
Because we defined $\eps_n(\cdot)$ in terms of spanning trees which
differ from the MST by a {\em non-random} proportion of edges, we
need a detour to handle expectations over events of asymptotically zero probability.
We defer the proof.
\begin{Lemma}
\label{L1}
(a) For any sequence $T^*_n$ of spanning trees, the sequence
$n^{-1} \len(T^*_n)$ is uniformly integrable.\\
(b) There exist spanning trees $T^{\prime \prime}_n$ such that
 \[  |T^{\prime \prime}_n \setminus T_n| \geq
a_n \]
where $a_n/n \to 1/2$.
\end{Lemma}
Now consider the spanning tree $T^*_n$ defined to be
$T^{\prime}_n$ if
$ n^{-1} |T^\prime_n \setminus T_n| \geq \sfrac{1}{10} q(\delta) $
and to be
$T^{\prime \prime}_n$ if not.
It follows from (\ref{con-1},\ref{con-2}) and Lemma \ref{L1} that
\[ n^{-1}|T^*_n \setminus T_n| \geq \sfrac{1}{10} q(\delta) \quad \mbox { (for
large $n$)}\]
\[ \limsup_n n^{-1} \Eb(\len(T^*_n) - \len(T_n)) \leq \sfrac{1}{9} r(\delta) . \]
Then from the definitions of $q(\delta), r(\delta)$ and
the assumption that $f(\cdot)$ is bounded it is easy to check
\begin{equation}
 q(\delta) \sim c \delta, \quad
r(\delta) \sim \sfrac{1}{2} \delta q(\delta)
\mbox{ as } \delta \downarrow 0 \label{qr-limits}
\end{equation}
for a certain $0<c<\infty$.
This establishes the upper bound (a) in Theorem \ref{T1}.

{\bf Proof of Lemma \ref{L1}.} Part (a) is automatic because,
writing $\sum_e$ for the sum over all edges of $\Cbold^2_m$, the
sequence $n^{-1} \sum_e \xi_e$ is uniformly integrable. For (b),
note that the cube $\Cbold^2_m$ with $2m(m-1)$ edges can be regarded
as a subgraph of the discrete torus $\Zbold^2_m$ with $2m^2$ edges.
Take a uniform random spanning tree $\widetilde{\TT}_n$ on
$\Zbold^2_m$, delete edges not in $\Cbold^2_m$ and add back boundary
edges to make some (non-uniform) random spanning tree $\TT_n$ on
$\Cbold^2_m$. By symmetry of the torus we have $\Pb(e \in
\widetilde{\TT}_n) = \frac{m^2 -1}{2m^2}$ for each edge $e$ of the
torus, and it follows that $\Pb(e \in \TT_n) = \frac{m^2 -1}{2m^2}$
for each non-boundary edge of the cube. Since there are $4(m-1)$
boundary edges and $2(m-1)(m-2)$ non-boundary edges, for any
spanning tree ${\bf t}$ we have
\[ \Eb | \TT_n \cap {\bf t}| \leq 4(m-1) + (n-1)(m^2 -1)/(2m^2)
= 4(n^{1/2} -1) + (n-1)^2/(2n). \]
So
\begin{eqnarray*}
\Eb|\TT_n \setminus {\bf t}| &=& (n-1) - \Eb|\TT_n \cap {\bf t}| \\
&\geq& a_n:= (n-1) - 4(n^{1/2}-1) - (n-1)^2/(2n) .
\end{eqnarray*}
So for any spanning tree ${\bf t}$ there exists some spanning tree
${\bf t}^*$ such that
$|{\bf t}^* \setminus {\bf t}| \geq a_n $.
Applying this fact to the MST gives (b).

\subsection{Upper bound: other cases}
The argument for Model 1 in the case $d \geq 3$ involves
only very minor modifications of the proof above, so we turn to
Model 2 with $d=2$ (the case $d \geq 3$ is similar).
Here it is natural to consider a different notion of
special configuration.

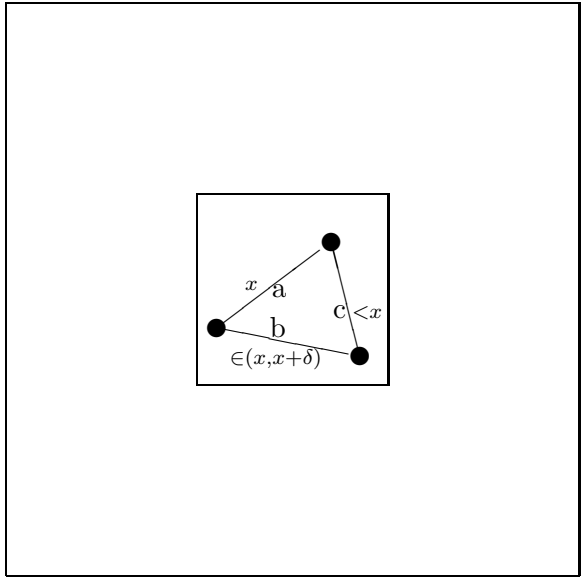
\begin{figure}
\begin{center}
\setlength{\unitlength}{0.1in}
\begin{picture}(30,33)(0,0)
\put(0,0){\line(1,0){30}}
\put(0,30){\line(1,0){30}}
\put(0,0){\line(0,1){30}}
\put(30,0){\line(0,1){30}}
\put(10,10){\line(1,0){10}}
\put(10,20){\line(1,0){10}}
\put(10,10){\line(0,1){10}}
\put(20,10){\line(0,1){10}}
\put(11,13){\circle*{1}}
\put(18.5,11.5){\circle*{1}}
\put(17,17.5){\circle*{1}}
\put(11,13){\line(5,-1){6.9}}
\put(11,13){\line(4,3){5.4}}
\put(18.5,11.5){\line(-1,4){1.5}}
\put(12.5,14.99){$\scriptstyle{x}$}
\put(11.7,11.1){$\scriptstyle{\in (x,x+\delta)}$}
\put(18.1,13.5){$\scriptstyle{<x}$}
\put(13.9,14.6){a}
\put(13.8,12.5){b}
\put(17.1,13.5){c}
\end{picture}
\end{center}
\caption{A special configuration on the $3 \times 3$ square.}
\end{figure}

\noindent
Here there is a $3 \times 3$ square containing a concentric
$1 \times 1$ square.
There are three points within the larger square, all being inside
the smaller square.
In the triangle $abc$ formed by the three points,
writing $x$ for the length of the second longest edge length,
the length of the longest edge is in the interval  $ (x, x+ \delta)$,  and $x+ \delta < 1$.
For such a configuration (within a configuration on a $m \times m$ square containing the
$3 \times 3$ square), edges
$ac$ are in the MST, and edge $b$ is not.
We can modify the minimal spanning tree by removing edge $a$
and adding edge $b$; 
this creates a new spanning tree whose extra length equals
$\len(b) - x$.

We now repeat the argument from the previous section, and the overall logic is
the same.
One gets different formulas for $q(\delta), r(\delta)$
but they have the same relationship (\ref{qr-limits}).
The weak law (\ref{con-1},\ref{con-2}) is easily established.
The only non-trivial difference is that we need to replace the
technical Lemma \ref{L1}
by the following technical lemma.
\begin{Lemma}
(a) There exists $c_1$ such that for any $n$ and any configuration
on $n$ points in the square of area $n$, the MST $\hat{T}_n$
has
$\len(\hat{T}_n) \leq c_1 n$.\\
(b) For sufficiently large $n$,
there exist spanning trees
$T^{\prime \prime}_n$ such that
$\len(T^{\prime \prime}_n) \leq
12c_1n$
and
 \[ n^{-1} |T^{\prime \prime}_n \setminus T_n| \geq
\sfrac{1}{2} . \]
\end{Lemma}
{\bf Proof.}
Part (a) follows from the analogous result for TSP -- see
\cite{steele97} inequality (2.14).
For (b), let $\xi_1,\ldots,\xi_n$
be the positions of the $n$ random points
and recall that $T_n$ is their MST.
Classify these points $\xi_i$ as ``odd" or ``even" according
to whether the number of edges in the path inside $T_n$
from $\xi_i$ to $\xi_1$ is odd or even.
Let $(\hat{\xi}_i)$ be a configuration obtained from $(\xi_i)$
by moving each ``odd" point a distance $11c_1$ in some arbitrary
direction.
Let $\hat{T}_n$ be the MST on $(\hat{\xi}_i)$.
Let $T^{\prime \prime}_n$ be the spanning tree on $(\xi_i)$ defined by
\[ (\xi_i,\xi_j) \in T^{\prime \prime}_n
\mbox{ iff }
(\hat{\xi}_i,\hat{\xi}_j) \in \hat{T}_n . \]
Suppose
$(\xi_i,\xi_j)$
is an edge of both $T_n$ and $T^{\prime \prime}_n$.
Since one end-vertex is odd and the other is even,
it is easy to see:\\
either (i) $\len(\xi_i,\xi_j) \geq 5 c_1$;
or (ii) $\len(\hat{\xi}_i,\hat{\xi}_j) \geq 5c_1$.\\
But by part (a) there are at most $n/5$ edges satisfying (i),
and similarly for (ii).
So $|T_n \cap T^{\prime \prime}_n| \leq 2n/5$.
Noting that
\[ \len(T^{\prime \prime}_n) \leq 11c_1 (n-1) + \len(\hat{T}_n) \leq 12c_1n \]
using (a), we have established (b).

\subsection{The lower bound: a discrete lemma}
The lower bound argument rests upon the following simple lemma.
\begin{Lemma}
\label{L2}
Consider a finite connected network with distinct edge-lengths.
If $T$ is the MST and $T^\prime$ is any ST then
\[ \len(T^\prime) - \len(T) \geq \sum_{e^\prime \in T^\prime \setminus T} \exc
(e^\prime) . \]
\end{Lemma}
{\bf Proof.}
Suppose
$ |T^\prime \setminus T| = k \geq 1$.
It is enough to show that there exist
$e^\prime \in T^\prime \setminus T$
and
$e \in T \setminus T^\prime$
such that\\
(i) $T^* = T^\prime \setminus \{e^\prime\} \cup \{e\}$
is a ST;\\
(ii) $|T^* \setminus T| = k-1$;\\
(iii) $\len(e^\prime) - \len(e) \geq \exc(e^\prime)$\\
for then we can continue inductively.

To prove this we first choose an arbitrary
$e^\prime \in T^\prime \setminus T$.
Consider $T^\prime \setminus \{e^\prime\}$.
This is a two-component forest; 
so the path in $T$ linking the end-vertices of $e^\prime$ must
contain some edge $e \in T \setminus T^\prime$ which links these two
components. So choose some such edge $e$. Properties (i) and (ii)
are clear. Apply Lemma \ref{L0} (b) to the end-vertices of
$e^\prime$ to see that $\perc(e^\prime) \geq \len(e)$. So
\[ \len(e^\prime) - \len(e) \geq \len(e^\prime) - \perc(e^\prime) = \exc(e^\prime) \]
which is (iii).

We will also need the following integration lemma;
part (a) will be used for Model 1 and part (b) for Model 2.
\begin{Lemma}
\label{LxiW}
(a)  Let $\xi$ and $W$ be independent real-valued random variables
such that $\xi$ has a density function bounded by a constant
$b$.
Then for any event $A \subseteq \{\xi > W\}$ we have
\[ \Eb(\xi - W)\ind_A \geq \sfrac{\Pb^2(A)}{2 b} . \]
(b)  Let $(V_i, i \geq 1)$ be real-valued r.v.'s such that
$\mu(0,x) := \Eb \sum_i \ind_{(0<V_i<x)}$ satisfies $\limsup_{x
\downarrow 0} \frac{\mu(0,x)}{x} < \infty$. Then there exists a
function $g(s) \sim \beta s^2$ as $s \downarrow 0$, for some $\beta
> 0$, such that for any sequence of events $A_i \subseteq \{V_i >
0\}$,
\[ \Eb \sum_i V_i\ind_{A_i} \geq
g \left( \sum_i \Pb(A_i) \right) . \]
\end{Lemma}
{\bf Proof.} (a)  It is sufficient to prove
\begin{eqnarray}
\label{ineq:cond} \Eb[(\xi - W)\ind_A|W]\geq \sfrac{\Pb^2(A|W)}{2
b}\mbox{ a.s.},
\end{eqnarray}
then by Jensen's inequality, we get
\[
\Eb[\Eb[(\xi - W)\ind_A|W]]\geq \sfrac{\Eb[\Pb^2(A|W)]}{2 b}\\
\geq \sfrac{\Pb^2(A)}{2 b}.
\]
Since $\xi$ and $W$ are independent of each other, equation
(\ref{ineq:cond}) reduces to
\[
\Eb[\xi \ind_A]\geq \sfrac{\Pb^2(A)}{2 b},
\mbox{ for } A\subseteq \{\xi>0\} .
\]
We can couple $\xi$ to a r.v. $U$ such that \\
(i) $U = 0$ on $\{ \xi \leq 0\}$.\\
(ii) $U$ has constant density $b$ on $(0, \Pb(\xi > 0)/b)$.\\
(iii) $U \leq \xi$.\\
Now it suffices to prove
\begin{equation}
\Eb[U \ind_A]\geq \sfrac{\Pb^2(A)}{2 b},
\mbox{ for } A\subseteq \{U>0\} . \label{UA}
\end{equation}
But it is clear that, for a given value of $\Pb(A)$, the choice of $A\subseteq \{U>0\} $ that minimizes
$\Eb [U\ind_A]$ is of the form
$A_c := \{0<U<c\}$ for some $c< \Pb(\xi>0)/b$.
A brief calculation gives
\[ \Eb[U \ind_{A_c}] =  \sfrac{\Pb^2(A_c)}{2 b} \]
establishing (\ref{UA}).

(b) For small $s > 0$ define $g(s)$ by
\[ g(s) = \int_0^{c(s)} x \mu(dx)
\mbox{ where }
\mu(0,c(s)) = s . \]
By hypothesis there exists $\gamma > 0 $ such that
$c(s) \geq \gamma s$ for small $s$.
So
\[ g(s) \geq \int_{c(s/2)}^{c(s)} x \mu(dx)
\geq c(s/2) \times s/2
= \gamma s^2/4 . \]
Taking
$A_i^s := \{0<V_i \leq c(s)\}$ we have
\[ \sum_i \Pb(A_i^s) = \mu(0,c(s)) = s; \quad
 \Eb \sum_i V_i\ind_{A^s_i} =
\int_0^{c(s)} x \mu(dx) = g(s) . \]
This is clearly the choice of $(A_i)$
which minmizes the left side subject to $\sum_i \Pb(A_i) = s$, and so for arbitrary $(A_i)$ we have
\[ \Eb \sum_i V_i\ind_{A_i} \geq
g \left( \sum_i \Pb(A_i) \right) . \]

\subsection{The lower bound in Model 1}
\label{sec-LB1}
We treat the case $d = 2$, but $d \geq 3$ involves only minor changes.
Recall $\Cbold_m^2 = G^{(n)}$ has
$c_n := 2(n - n^{1/2})$ edges.
Fix $\delta > 0$.
Consider a pair
$(T^\prime_n,T_n)$ attaining the minimum in
the definition (\ref{eps-def}) of $\eps_n(\delta)$.
For a uniform random edge $e_n$ of $\Cbold_m^2$,
\begin{equation}
\Pb(e_n \in T^\prime_n \setminus T_n) = \frac{\Eb|T^\prime_n
\setminus T_n|}{c_n} \geq \frac{\delta}{2} \label{ett}
\end{equation}
and
\begin{eqnarray}
\Eb \eps_n(\delta)
&=& \sfrac{1}{n} \Eb(\len(T^\prime_n) - \len(T_n))\nonumber\\
&\geq& \sfrac{1}{n} \Eb \sum_{e \in T^\prime_n \setminus T_n} \exc
(e)
\mbox{ by Lemma \ref{L2}}\nonumber\\
&=&\frac{c_n}{n} \ \Eb \exc(e_n)1_{(e_n \in T^\prime_n \setminus
T_n)} \label{eee}
\end{eqnarray}
For a fixed edge $e$ of $\Cbold_m^2$
we can write
\[ \exc(e) = (\xi^{(n)}(e) - W^{(n)}(e))^+ \]
where $\xi^{(n)}(e)$ is the edge-length of $e = (v,v^*)$
and where $$W^{(n)}(e) =  \inf \{ t:
\mbox{ $v$ and $v^*$  in the same component of $G^{(n)}_t \setminus \{ e \} $}\}.$$
Note (and this is the key special feature that makes Model 1
easy to study)
that $\xi^{(n)}(e)$ and $W^{(n)}(e)$ are independent.
Since $\exc(e) > 0 $ on $\{e \in T^\prime_n \setminus T_n\}$
we see that the quantity at (\ref{eee}) is of the form
appearing in Lemma \ref{LxiW}(a).
So
\begin{eqnarray*}
\frac{n}{c_n} \Eb \eps_n(\delta) &\geq & \Eb\left(\xi^{(n)}(e_n) -
W^{(n)}(e_n)\right)1_{(e_n \in T^\prime_n \setminus T_n)}
\mbox{ by (\ref{eee})}\\
&\geq & \frac{\Pb^2(e_n \in T^\prime_n \setminus T_n)}{2 \bar{f}}
\mbox{ by Lemma \ref{LxiW}(a)}\\
&\geq& \frac{\delta^2}{8 \bar{f}} \mbox{ by (\ref{ett})}
\end{eqnarray*}
where $\bar{f}$ is the bound on the density of $\xi$.
Because $c_n \sim 2n$
we have established part (b) of Theorem \ref{T1} in this case.

\section{The minimum spanning forest and continuum percolation}\label{sec:infinite}
It remains to prove the lower bound in Model 2.
Rather than doing calculations with the finite model,
we consider the limit Poisson process on the plane, and exploit the well known connection
between the
{\em minimum
spanning forest} (MSF) 
and continuum percolation.
We then relate the finite models to the infinite limits in section \ref{sec-LB2},
as an instance of {\em local weak convergence} \cite{me101} of random graphical structures.

\subsection{Minimum spanning forests}\label{sec:msf}

Here is a general definition, in the context of a countable-vertex
network $G$ with distinct edge-lengths (see \cite{alexander95} for
more detailed treatment). As in Section \ref{sec-back} let $G_t$ be
the subnetwork consisting of those edges $e$ of $G$ with $\len (e) <
t$. Define the MSF by:
\begin{quote}
an edge $(v,w)$ is in the MSF if and only if, for
$t = \len(v,w)$, vertices $v$ and $w$ are in different components
of $G_t$ and at least one of these components is finite.
\end{quote}

Consider a
Poisson point process $\Phi = \sum_i \delta_{\eta_i}$ of rate $1$ in
$\Rbold^d$. Add an extra point $O$ at the origin. Consider $\Phi^O =
\sum_i \delta_{\eta_i} +\delta_{O}$ as the vertices of a network G
(the complete graph with Euclidean edge-lengths). With probability
one, $\Phi^O$ has only finitely many points in any bounded subset of
$\Rbold^d$ and all of the interpoint distances are distinct. As in
Section \ref{sec-back}, we define for arbitrary points $\eta_i$ and
$\eta_j$ of $\Phi^O$,
\begin{eqnarray*}
\perc(\eta_i,\eta_j) =\inf\{t: \mbox{ $\eta_i$ and $\eta_j$ are in
the same component of $G_t$}\}.
\end{eqnarray*}

We now give some properties of the MSF denoted $\mathcal{F}_\infty$ on this network and show how
Lemma \ref{L0} extends to this setting.

\begin{Lemma}\label{lem:percinfty}
(a) We have
$e\in \mathcal{F}_\infty$ if and only if $\len(e) = \perc(e)$.\\
(b)
For any vertex-pair $u,v$ write $u\to v$ for the set of paths $\pi$
from $u$ to $v$.  Then, a.s. 
\begin{eqnarray}
\label{perc:minimax} \perc(u,v)= \min_{\pi: u\to v}\max\{\len(e) : e
\in \pi\}
\end{eqnarray}
\end{Lemma}
{\bf Proof.} 
Let us say that $G$ has the
uniqueness property if for every vertex-pair $u,v\in G$, the graph
$G_{\len(u,v)}$ has at most one infinite component (note that this
notion was used in the proof of Lemma 2.1 in
\cite{penroseyukich}). Part (a) will follow from the fact that
$\Phi^O$ has the uniqueness property, which implies:
\begin{eqnarray*}
e=(v,w) \in \mathcal{F}_\infty &\Leftrightarrow & \mbox{$v$ and $w$
are in different components of $G_{\len(e)}$}\\
&\Leftrightarrow & \perc(e) \geq \len(e)\\
&\Leftrightarrow & \perc(e) = \len(e).
\end{eqnarray*}
To show that $\Phi^O$ has the uniqueness property almost surely
it is enough to show
\begin{eqnarray}
\label{eq:uniq}\Pb(\forall u\in \Phi^O, G_{\len(O,u)} \mbox{ has at
most one infinite component}) = 1.
\end{eqnarray}
This last fact follows from Theorem 1.8 (and Remark 1.10) of
\cite{alexander-uniq}, which implies (see also \cite{alexander95}),
\begin{eqnarray}
\label{eq:uniq-a}\Pb(G_t \mbox{ includes at most one infinite component for each } t\in
\Rbold) =1.
\end{eqnarray}
Note that (\ref{eq:uniq}) can be proved without appealing to the simultaneous
uniqueness result as follows:
\begin{eqnarray*}
\lefteqn{\Pb(\forall u\in \Phi^O, G_{\len(O,u)} \mbox{ has at most
one
infinite component})}\\
&=& \lim_{n\to \infty}\Pb(\forall u\in \Phi^O\cap
B(n), G_{\len(O,u)} \mbox{ has at most one infinite component})\\
&\geq& \lim_{n\to \infty}\Pb(\forall u\in \Phi^O\cap B(n),
G_{\len(O,u)}\setminus B(n) \mbox{ has at most one infinite
component}),
\end{eqnarray*}
where $B(n)$ is the ball of center the origin and radius $n$ and for
any network $G$ on $\Rbold^d$, $G\setminus B(n)$ is the subnetwork
with edges and vertices in $\Rbold^d\setminus B(n)$. By independence
and the fact that there can be at most one infinite component in
continuum percolation (see Theorem 3.6 in \cite{meesterroy}), we
have
\begin{eqnarray*}
\Pb(\forall u\in \Phi^O\cap B(n), G_{\len(O,u)}\setminus B(n) \mbox{
has at most one infinite component}|\Phi^O\cap B(n)) = 1 \mbox{
a.s.}
\end{eqnarray*}
which proves (\ref{eq:uniq}). 

We now prove (b). Let $t = \perc(u,v)$, the definition of $\perc(u,v)$ may be restated easily as:
\begin{eqnarray*}
t = \perc(u,v)= \inf_{\pi: u\to v}\max\{\len(e) : e
\in \pi\}
\end{eqnarray*}
Hence (b) amounts to prove that with probability one, this infimum is indeed a minimum. Note that $G_t(u)\cap G_t(v) =\emptyset$ and by
(\ref{eq:uniq-a}) a.s. at least one of these two clusters, say $G_t(u)$, is finite.
Let $E$ be the set of edges with exactly one of its end vertices in $G_t(u)$ and the other one in $G_t(u)^c$, 
and with edge length less than $t+1$. The
set $E$ is a.s. finite and then we easily see that $\min\{\len(e), e\in E\}
=t=\perc(u,v)$ since $u$ and $v$ are in the same component of
$G_{t+\epsilon}$ for any $\epsilon>0$. 
Let $e^*=\arg\min\{\len(e), e\in E\}$ and write $e^*=(a,b)$ with $a\in
G_t(u)$ and $b \in G_t(u)^c$. Since a.s. we have $\len(e)\neq \len(e^*)=t$ for any $e\neq e^*$, a.s. we
have $G_{t+}(u):=\cap_{\epsilon>0}
G_{t+\epsilon}(u)=G_t(u)\cup\{e^*\}\cup G_t(b)$
and $G_{t+}(v)=G_t(v)$. But the definition of $t=\perc(u,v)$ implies that
$G_{t+}(u) = G_{t+}(v)$, and hence $b\in G_t(v)$.
It follows that
\begin{eqnarray*}
\inf_{\pi:  u\to v}\max\{\len(e) : e \in \pi\}= \len(e^*) = \min_{\pi:  u\to v}\max\{\len(e) : e \in \pi\},
\end{eqnarray*}
and (b) follows.

\subsection{Finite density}

Define the measure $\mu$ on $(0,+\infty)$ by
\begin{eqnarray*}
\mu(0,x) &=& \Eb \sum_i \ind(0<\len(O,\eta_i)-\perc(O,\eta_i)<x)\\
&=& \Eb \sum_i \ind(0<\len(O,\eta_i)-\perc(O,\eta_i)\leq x).
\end{eqnarray*}
The next lemma formalizes the heuristic idea
$f_{\mu}(0^+) < \infty$ from section \ref{sec-heuristic}. 
\begin{Proposition}
\label{A1}
In Model 2, we have,
\[ \limsup_{x \downarrow 0}  \frac{\mu(0,x)}{x} \ < \infty . \]
\end{Proposition}

For $(X_1,\cdots, X_n) \in (\Rbold^d)^n$ we define
$\Phi^{X_1,\cdots, X_n} = \Phi + \sum_{i=1} ^n \delta_{X_{i}}$, and write
$\Pb^{X_1,\cdots, X_n}$ for the probability measure associated with the
random variable $\Phi^{X_1,\cdots, X_n}$. Using Campbell's
formula, we have
\begin{eqnarray*}
\mu(0,x) &=& \omega_d \int_0^\infty \Pb^{O,\ut}(\perc(O,\ut) \in
[t-x,t))t^{d-1} dt,
\end{eqnarray*}
where $\ut $ is the point $ (t,0,\dots,0)$ and
$\omega_d=\frac{2\pi^{d/2}}{\Gamma(d/2)}$ is the surface of the unit
sphere.

We need to introduce some continuum percolation terminology.  For any
$r$ and $\lambda$, we define the probability measure
$\Pb^{O,\ut}_{r}$ under which $\Phi$ is a Poisson point
process of intensity $1$ and an edge $e$ from the complete
graph $\Phi^{O,\ut}$ is said to be open (resp. closed) if
$\len(e)<r$ (resp. $\len(e)\geq r$). We denote by $G^O$ the open
cluster containing the origin: $G^O=G^{O}_r$. Let $r_c$ be the critical radius for the Poisson
continuum percolation model of density $1$ and deterministic
radius, i.e. for $r<r_c$ the number of vertices in any open
cluster is finite whereas for $r>r_c$ there exists an
unique unbounded open cluster.

Write $C, C_1, C_2$ for positive
constants not depending on the parameters of the problem.

\begin{Lemma}
\label{lem:outrc} For any $\epsilon>0$, we have
\begin{eqnarray*}
\mbox{for $0<t<r_c-\epsilon$, }&& \Pb^{O,\ut}(\perc(O,\ut)\in
[t-x,t))\leq C_1x,\\
\mbox{for $t>r_c+\epsilon$, }&& \Pb^{O,\ut}(\perc(O,\ut)\in
[t-x,t))\leq C_1xe^{-C_2 t}.
\end{eqnarray*}
\end{Lemma}
We first introduce some notations. The edge-length is the Euclidean distance
denoted $\len(u,v)=|u-v|$. For a set $S\subset \Rbold^d$, we denote
by $d(S) = \sup\{|x-y|,\:x,y\in S \}$ its diameter. For $x\in
\Rbold^d$ and $r>0$, $B(x,r)$ denotes the open ball of radius $r$
centered at $x$. For $t>0$ we denote $S(t)=[-t,t]^d$.
Under the probability measure $\Pb^{O,\ut}_{r}$, the
occupied region is $\cup_{X \in \Phi^{O,\ut}}B(X,r/2)$ and the
vacant region is the complement of the occupied region. The occupied
component of the origin $W$ is defined by
$\Pb^{O,\ut}_{r}(W=\cup_{X\in G^O}B(X,r/2))=1$. The vacant
component containing the point $\ut/2$ is denoted by $V$. More
generally, for $r>0$ the occupied region at level $r$ is $\cup_{X
\in \Phi^{O,\ut}}B(X,r/2)$ and we denote $W_r=\cup_{X\in
G^O_r}B(X,r/2)$ the occupied component of the origin at level $r$
and $V_r$ the vacant component containing the point $\ut/2$ at level
$r$.

Since we may assume that all interdistances are different, there
exists an unique pair $(X,Y)$ in the support of $\Phi^{O,\ut}$ such that $\perc(O,\ut)
= |X-Y|$ (see Lemma \ref{lem:percinfty}).

First consider the case $t<r_c-\epsilon$. Let $S_z= zt/2 + S(t/2)$
where $z\in \mathbb{Z}^d$. If the event $\{\perc(O,\ut)\in
[t-x,t)\}$ occurs, there is some $z\in \mathbb{Z}^d$ such that
$W\cap S_z\neq \emptyset$ and there exists $X,Y\in \Phi\cap S_z$
such that $|X-Y|\in [t-x,t)$. Note that we have for any $z\in
\mathbb{Z}^d$,
\begin{eqnarray*}
\Pb^{O,\ut}_{t}\left(\exists X,Y\in \Phi\cap S_z,\: |X-Y|\in
[t-x,t) \right) \leq C(1+t^{2d-1}) x.
\end{eqnarray*}
Hence we have
\begin{eqnarray}
\nonumber\lefteqn{\Pb^{O,\ut}(\perc(O,\ut)\in [t-x,t))}\\
\nonumber &\leq& \sum_{z\in
\mathbb{Z}^2}\Pb^{O,\ut}_{t}\left(W\cap S_z\neq \emptyset,\:
\exists X,Y\in \Phi(S_z),\: |X-Y|\in [t-x,t) \right)\\
\label{sum}&\leq& \left(K +\sum_{z,\:\|z\|\geq
2}\Pb^{O,\ut}_{t}\left(W\cap S_z\neq
\emptyset\right)\right)C(1+t^{2d-1}) x,
\end{eqnarray}
where  $\|z\|:=\max(|z_1|,|z_2|)$ and $K$ is a constant depending on
$d$. 
Lemma 3.3 of
\cite{meesterroy} ensures that the sum of (\ref{sum}) is finite for
$t<r_c-\epsilon$.

The case $t>r_c+\epsilon$ is quite similar. If the event
$\{\perc(O,\ut)\in [t-x,t)\}$ occurs, there is some $z\in
\mathbb{Z}^2$ such that $V_{t-x}\cap S_z\neq \emptyset$ and there
exists $X,Y\in \Phi\cap S_z$ such that $|X-Y|\in [t-x,t)$. Hence we
have
\begin{eqnarray}
\nonumber\Pb^{O,\ut}(\perc(O,\ut)\in [t-x,t))&\leq& \left(5
+\sum_{z,\:\|z-\ut/2\|\geq 2}\Pb^{O,\ut}_{t}\left(V_{t-x}\cap
S_z\neq
\emptyset\right)\right)C(1+t^{2d-1})x\\
\label{exp}&\leq& C_1e^{-C_2 t}x,
\end{eqnarray}
where (\ref{exp}) follows from Lemma 4.1 of \cite{meesterroy} and
the fact that
\begin{eqnarray*}
\Pb^{O,\ut}_{t}\left(V_{t-x}\cap S_z\neq \emptyset\right) \leq
\Pb^{O,\ut}_{t}\left(d(V_{t-x})>\|z-\ut/2\|\right).
\end{eqnarray*}

We now concentrate on the case $t\in (r_c-\epsilon,r_c+\epsilon)$.
We define the event $$A=\{ \mbox{the points of } \Phi \mbox{ on the
  axis $e_1$ are in } G^O\}, \mbox{ where $e_1=(1,0,\dots,0)$.}$$ Under $\Pb^{O,\ut}_{r}$, with probability one, we have
$A =\{ \ut \in G^O\}$ and
\begin{eqnarray}
\nonumber \Pb^{O,\ut}(\perc(O,\ut)\in [t-x,t)) &=&
\Pb^{O,\ut}_{t} (A )-\Pb^{O,\ut}_{t-x}(A)\\
\label{eq:perc10}& = &  \Pb^{O,\ut}_{t} (A) -  \Pb^{O,\utx}_{t-x} (A) + \Pb^{O,\utx}_{t-x} (A)  -  \Pb^{O,\ut}_{t-x} (A).
\end{eqnarray}
We first prove that
\begin{equation}
\label{eq:ineq2}
\Bigm|  \Pb^{O,\utx}_{t-x} (A)  -  \Pb^{O,\ut}_{t-x} (A) \Bigm| \leq C t ^d x.
\end{equation}
Note that
\begin{eqnarray*}
\Pb^{O,\ut}_{t-x} ( A ) &=& \Pb^{O}_{t-x}
(B(\ut,t-x)\cap G^O\neq \emptyset )\\
\Pb^{O,\utx}_{t-x} ( A ) &=& \Pb^{O}_{t-x}
(B(\utx,t-x)\cap G^O\neq \emptyset ),
\end{eqnarray*}
where $B(X,r)$ denotes the open ball of radius $r>0$ centered at
$X\in \Rbold^d$. Hence we have
\begin{eqnarray*}
\left| \Pb^{O,\utx}_{t-x} ( A ) - \Pb^{O,
\ut}_{t-x} (A)\right| \leq \Pb(\Phi(B(\ut,t-x) \Delta
B(\utx, t-x))\geq 1)\leq Ct^d x,
\end{eqnarray*}
where $B(\ut,t) \Delta B( \utx, t-x)$
denotes the symmetric difference.  This is exactly (\ref{eq:ineq2}).

We then write:
\begin{eqnarray}
\frac{1}{x} \int_{r_c - \epsilon} ^ {r_c + \epsilon} \Pb^{O,\ut} (\perc(O,\ut) \in [t-x,t))  t^{d-1} dt &  =  & 
\frac{1}{x} \int_{r_c - \epsilon} ^ {r_c + \epsilon} \Pb^{O,\ut}_{t} (A) t^{d-1} dt 
-   \frac{1}{x} \int_{r_c - \epsilon} ^ {r_c + \epsilon}  \Pb^{O,\utx}_{t-x} (A)t^{d-1} dt \nonumber \\
 &&  \quad+ \frac{1}{x} \int_{r_c - \epsilon} ^ {r_c + \epsilon}  (  \Pb^{O,\utx}_{t-x} (A)  
 -  \Pb^{O,\ut}_{t-x} (A) )t^{d-1} dt  \label{eq:decomp}
\end{eqnarray}

With the change of variable $t \mapsto t- x$, the second term on the right hand side  of (\ref{eq:decomp}) is decomposed as follows:
$$
 \frac{1}{x} \int_{r_c - \epsilon} ^ {r_c + \epsilon}  \Pb^{O,\utx}_{t-x} (A)t^{d-1} dt  \leq  
 \frac{1}{x} \int_{r_c - \epsilon -x} ^ {r_c + \epsilon-x}  \Pb^{O,\ut}_{t} (A)t^{d-1} dt  
 +  K\int_{r_c - \epsilon -x} ^ {r_c + \epsilon-x}  \Pb^{O,\ut}_{t} (A) t^{d-2}dt, 
$$
where $K$ is a constant depending on $d$.

Hence,  the decomposition (\ref{eq:decomp}) is further decomposed as
\begin{eqnarray*}
\Bigm|\frac{1}{x} \int_{r_c - \epsilon} ^ {r_c + \epsilon} \Pb^{O,\ut} (\perc(O,\ut) \in [t-x,t))  t^{d-1} dt \Bigm|
&  \leq  &  \frac{1}{x} \int_{r_c + \epsilon-x} ^ {r_c + \epsilon} \Pb^{O,\ut}_{t} (A) t^{d-1} dt  
+  \frac{1}{x} \int_{r_c - \epsilon-x} ^ {r_c - \epsilon} \Pb^{O,\ut}_{t} (A) t^{d-1} dt \\
&& + K\int_{r_c - \epsilon -x} ^ {r_c + \epsilon-x}  \Pb^{O,\ut}_{t} (A)t^{d-2} dt \\
&& + \frac{1}{x} \int_{r_c - \epsilon} ^ {r_c + \epsilon}  (  \Pb^{O,\utx}_{t-x} (A)  
-  \Pb^{O,\ut}_{t-x} (A) )t^{d-1} dt
\end{eqnarray*}
By (\ref{eq:ineq2}), the last term is bounded by $ \int_{r_c - \epsilon} ^ {r_c + \epsilon} Ct^{2d-1}dt = C_1 $. It implies that 
\begin{eqnarray}
\label{eq:rc}\Bigm| \frac{1}{x} \int_{r_c - \epsilon} ^ {r_c + \epsilon} \Pb^{O,\ut} (\perc(O,\ut) \in [t-x,t))  t^{d-1} dt \Bigm| \leq C_2.
\end{eqnarray}

Proposition \ref{A1} now follows from Lemma \ref{lem:outrc} and Equation
(\ref{eq:rc}).

\subsection{The lower bound in Model 2}
\label{sec-LB2}

We start with a slight extension of Proposition 9 of \cite{me57}
(see also Theorem 7 in \cite{me101}). In what follows, a set of
points is identified with its associated geometric graph which is
the complete graph over these points with Euclidean distance as
edge-lengths.
\begin{Lemma}
Let $\Phi_n$ denote the point process consisting of $n$ points
$\{\xi_i, \: 1\leq i\leq n\}$ which are independent and have the
uniform distribution on the square $[0,n^{1/d}]^d$. For each $n$,
let $U_n$ be chosen independently and uniformly from the set
$\{1,\dots, n\}$, and let
\begin{eqnarray*}
\Phi^O_n = \{ \xi^{(n)}_i := \xi_i-\xi_{U_n},\: 1\leq i \leq n\}.
\end{eqnarray*}
To each vertex $\xi^{(n)}_i$ of the rooted (at the origin) geometric
graph $\Phi^O_n$, we associate the mark
$\perc^{n}_i=\perc(O,\xi^{(n)}_i)$ as defined in (\ref{def-perc}).
We denote by $(\Phi^O_n, \perc^n)=\{\xi^{(n)}_i, \perc^{n}_i\}$ the
corresponding marked geometric graph.  Then one has joint weak
convergence
\begin{eqnarray}\label{lem:conv}
( (\Phi^O_n, \perc^n), MST(\Phi^O_n))\cd ((\Phi^O,\perc) ,
\mathcal{F}_\infty),
\end{eqnarray}
where $(\Phi^O,\perc)$ is the Palm version of the Poisson process
of intensity $1$ with the mark $\perc(O,\eta_i)$ associated to point
$\eta_i$.
\end{Lemma}
Here convergence
$MST(\Phi^O_n) \cd \mathcal{F}_\infty$
is {\em local weak convergence} in the sense of \cite{me101}.

{\bf Proof.}  The analog of (\ref{lem:conv}) without marks is Proposition 9 of
\cite{me57}. 
By the Skorokhod representation theorem, we can assume that with
probability one, we have
\begin{eqnarray}\label{con-as}
( \Phi^O_n=\{\xi^{(n)}_i\}, MST(\Phi^O_n))\rightarrow (\Phi^O =
\{\eta_i\}, \mathcal{F}_\infty)
\end{eqnarray}
We have to prove that for any $i\geq 1$,
\begin{eqnarray*}
\lim_{n\to \infty} \perc(O,\xi^{(n)}_i)=\perc(O,\eta_i) \mbox{ a.s.}
\end{eqnarray*}
By Lemma \ref{lem:percinfty}, we know that
$\perc(O,\eta_i)=\max\{\len(e),\: e\in \pi^*\}$ where $\pi^*$ is the minimax
path from $O$ to $\eta_i$. By definition of
the metric of local weak convergence, (\ref{con-as}) implies
that for arbitrary fixed $L$, we have with $S(L)=[-L,L]^d$,
\begin{eqnarray*}
\forall \eta_i\in S(L),\quad \xi^{(n)}_i \to \eta_i.
\end{eqnarray*}
For $L$ sufficiently large, the path $\pi^*$ is included in $S(L)$
and let $\pi^*_n$ be the associated path in $\Phi^O_n$ Since
$$\perc(O,\xi^{(n)}_i) =\min_{\pi:O\to \xi^{(n)}_i}\max_{e\in
\pi}\len(e)\leq \max_{e\in \pi^*_n}\len(e),$$ by the convergence of
$\pi^*_n$ to $\pi^*$, we have
\begin{eqnarray*}
\limsup_{n\to \infty} \perc(O,\xi^{(n)}_i) \leq \max\{\len(e),\:
e\in \pi^*\}=\perc(O,\eta_i).
\end{eqnarray*}
Now we need to prove that
\begin{eqnarray}
\label{per-liminf}\liminf_{n\to \infty} \perc(O,\xi^{(n)}_i) \geq
\perc(O,\eta_i).
\end{eqnarray}
Take $\pi^{(n)}:O\to \xi^{(n)}_i$ such that
\begin{eqnarray*}
\max _{e\in \pi^{(n)}} \len(e) =\perc(O,\xi^{(n)}_i).
\end{eqnarray*}
For $r>0$, we denote by $G_r(\eta_i)$ (resp. $G^n_r(\xi_i^{(n)})$)
the connected component of $\Phi^O$ (resp. $\Phi^O_n$) with edge
length less than $r$ containing $\eta_i$ (resp. $\xi_i^{(n)}$). Let
$\perc(O,\eta_i)=t$, so that we have $G_t(O)\cap
G_t(\eta_i)=\emptyset$ and say $G_t(O)$ is finite (see the
uniqueness property in the proof of Lemma \ref{lem:percinfty}). 
We define $\tilde{G}_{t}$ (resp.
$\tilde{G}^n_{t}$) to be the subgraph of $\Phi^O$ (resp.
$\Phi^O_n$) consisting of those edges with length less than
$t+1$ with exactly one of its end vertices in $G_t(O)$ (resp.
$G^n_t(O)$).
Let $e^*=\arg\max\{ \len(e), e\in \pi^*\}$. By Lemma \ref{lem:percinfty}, we know that
$\perc(O,\eta_i)=\len(e^*)$ and $e^*\in\tilde{G}_{t}$ is such that $e^* = \arg\min\{\len(e),\: e\in
\tilde{G}_{t}\}$. Since $G_t(O)$ is finite, we have clearly
that $\tilde{G}_{t}$ is included in $S(L)$ for sufficiently
large $L$. Then we have
\begin{eqnarray*}
\max _{e\in \pi^{(n)}} \len(e)&\geq& \min \{ \len(e),\: e\in
\tilde{G}^n_{t}\}\to t=\perc(O,\eta_i) \mbox{ as
$n\to\infty$},
\end{eqnarray*}
where the last limit follows from the convergence of
$\tilde{G}^n_{t}$ to $\tilde{G}_{t}$.

We now return to the proof of the lower bound in Model 2. We start
by copying and modifying the argument from section \ref{sec-LB1}.
Fix $\delta > 0$. Let $\xi_{U_n}$ be a uniform random vertex from
$(\xi_i, 1 \leq i \leq n)$. Consider a pair $(T^\prime_n,T_n)$
attaining the minimum in the definition (\ref{eps-def}) of
$\eps_n(\delta)$. Then
\begin{equation}
\Eb \sum_i \ind\{
 (\xi_{U_n},\xi_i) \in T^\prime_n \setminus T_n\}
 = \frac{2\Eb|T^\prime_n \setminus T_n|}{n} \geq 2 \delta
\label{xiU}
\end{equation}
and
\begin{eqnarray}
\Eb \eps_n(\delta)
&=& \sfrac{1}{n}\Eb(\len(T^\prime_n) - \len(T_n))\nonumber\\
&\geq& \sfrac{1}{n}\Eb \sum_{e \in T^\prime_n \setminus T_n} \exc
(e)
\mbox{ by Lemma \ref{L2}}\nonumber\\
&=&\sfrac{1}{2} \Eb \sum_i
\exc(\xi_{U_n},\xi_i)\ind \{(\xi_{U_n},\xi_i) \in T^\prime_n \setminus
T_n\} \label{eee2}
\end{eqnarray}
Note
that for $0 < L < \infty$
\begin{eqnarray*}
\Eb \sum_i \ind \{\len(\xi_{U_n},\xi_i) \geq L, (\xi_{U_n},\xi_i) \in
T_n^\prime\}
&=& \frac{2}{n} \Eb |\{e \in T^\prime_n: \len(e) \geq L\}|\\
&\leq & \frac{2}{n} \ \frac{\Eb \  \len(T^\prime_n)}{L} \\
&\leq & \delta
\mbox{ for $L = L(\delta)$ sufficiently large}
\end{eqnarray*}
the last inequality because $\Eb \len(T^\prime_n) = O(n)$. So fixing
such an $L$, (\ref{xiU}) implies
\begin{equation}
\Eb \sum_i \ind\{
 (\xi_{U_n},\xi_i) \in T^\prime_n \setminus T_n
, \len(\xi_{U_n},\xi_i) \leq L \}
\geq  \delta
\label{eLd}
\end{equation}
while (\ref{eee2}) trivially implies
\begin{equation}
2 \Eb \eps_n(\delta) \geq \ \Eb \sum_i
\exc(\xi_{U_n},\xi_i)\ind \{(\xi_{U_n},\xi_i) \in T^\prime_n \setminus
T_n, , \len(\xi_{U_n},\xi_i) \leq L\} . \label{eee3}
\end{equation}
The purpose of these representations is to exploit local weak
convergence.
Consider the near-minimal STs $T^\prime_n$ appearing in
(\ref{eLd},\ref{eee3}).  By a compactness argument and by passing to
a subsequence of $n$ we may assume that they converge to some forest
$\FF^\prime_\infty$ on $(\bar{\eta}_i)$; that is, we may assume that
(\ref{lem:conv}) remains true when we append $T^\prime_n $ to the
left side and $ \FF^\prime_\infty $ to the right side. We can now
take limits in (\ref{eLd}) to deduce
\[
\sum_i \Pb( (O,\eta_i) \in \FF^\prime_\infty \setminus
\FF_\infty, \ \len (O,\eta_i) \leq L) \geq \delta .
\] And taking limits in (\ref{eee3}) gives
\begin{equation} 2 \liminf_n \Eb \eps_n(\delta) \geq \Eb
\sum_i (\len (O,\eta_i) - \perc (O,\eta_i) )\ind \{
(O,\eta_i) \in \FF^\prime_\infty \setminus \FF_\infty, \
\len (O,\eta_i) \leq L\} . \label{qwe}
\end{equation}
Writing
\[ V_i =
 \len
(O,\eta_i) - \perc (O,\eta_i)
\] \[ A_i = \{
(O,\eta_i) \in \FF^\prime_\infty \setminus \FF_\infty, \
\len (O,\eta_i) \leq L)\} , \] we are precisely in the
setting in which Proposition \ref{A1} and Lemma \ref{LxiW}(b) apply, and the
conclusion is that the right side of (\ref{qwe}) is $\geq (\beta -
o(1)) \delta^2$ for small $\delta$, implying the lower bound in
Theorem \ref{T1}.


\end{document}